\documentstyle[12pt]{article}
\def\Bbb R{{\rm \bf R}}
\def\proclaim#1{\vskip2mm{\bf #1}\em}
\def\endproclaim{\em \vskip2mm}
\def\tag#1{\eqno(#1)}
\def\gathered{\begin{array}{c}}
\def\endgathered{\end{array}}
\def\text{\mbox}
\begin{document}

\title {Detecting a hidden obstacle via the time domain enclosure method. A scalar wave case}
\author{Masaru IKEHATA\footnote{
Laboratory of Mathematics,
Graduate School of Engineering,
Hiroshima University, Higashihiroshima 739-8527, JAPAN}}
\maketitle
\begin{abstract}
The characterization problem of the existence
of an unknown obstacle {\it behind} a known obstacle is considered
by using a singe observed wave at a place where the wave is generated.
The unknown obstacle is invisible from the place by using visible ray.
A mathematical formulation of the problem using the classical wave equation is given.
The main result consists of two parts:

(i)  one can make a decision whether the 
unknown obstacle exists or not behind a known impenetrable obstacle
by using a single wave over a finite time interval under some a-priori
information on the position of the unknown obstacle;

(ii) one can obtain a lower bound of the {\it Euclidean} distance of the unknown obstacle
to the center point of the support of the initial data of the wave.

The proof is based on the idea of the {\it time domain enclosure method}
and employs some previous results on the Gaussian lower/upper estimates for the {\it heat kernels}
and domination of semigroups.

\noindent
AMS: 35R30, 35L05, 35K08


\end{abstract}


\section{Introduction}
The purpose of this paper is to develop a {\it mathematical method} of imaging an unknow obstacle behind a known {\it impenetrable}
obstacle form a {\it single} wave generated and observed at the same place where one can not see the unknown obstacle
by using {\it visible ray}.  The governing equation of the wave should be the acoustic wave equation or the Maxwell system
for application.  In this paper, as a first step, we choose the classical wave equation.

Let us formulate our problem more precisely.
First of all we prepare three open sets described below.
Let $D_0$ be a nonempty bounded open subset of $\Bbb R^3$ with $C^2$-boundary such that $\Bbb R^3\setminus\overline{D_0}$ is connected.
Let  $D$ be a possibly empty bounded open subset of $\Bbb R^3$ with $C^2$-boundary such that $\overline{D_0}\cap\overline D=\emptyset$ and
$\Bbb R^3\setminus\overline D$ is connected.
Let $B$ be an open ball such that $\overline B\cap(\overline{D_0}\cup\overline D)=\emptyset$.
We think that the radius of $B$ is very small.

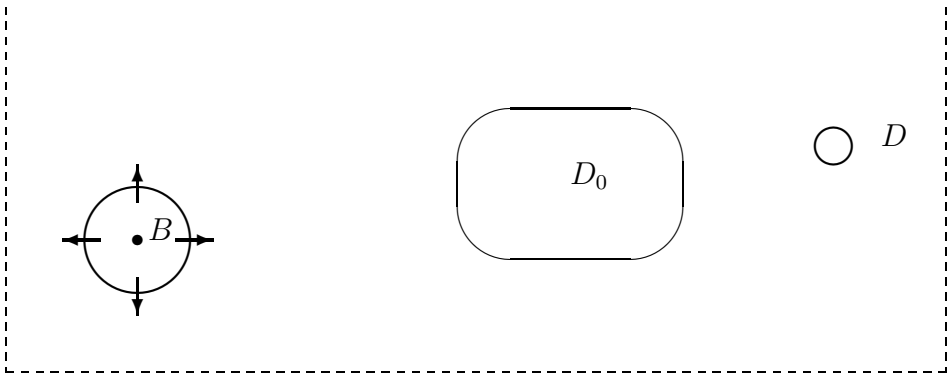
\begin{figure}
\setlength\unitlength{0.5truecm}
  \begin{picture}(10,10)(-3,-1)
  \put(0,0){\dashbox{0.2}(25,10){}}
  \put(15,5){\oval(6,4){$D_0$}}
  \put(3.5,3.5){\circle*{0.3}{$B$}}
  \thicklines
  \put(22,6){\circle{1} {$D$}}
  \put(3.5,3.5){\circle{3}}
  \put(4.5,3.5){\vector(1,0){1}}
  \put(3.5,4.5){\vector(0,1){1}}
  \put(3.5,2.5){\vector(0,-1){1}}
  \put(2.5,3.5){\vector(-1,0){1}}
  \end{picture}
 \caption{\label{fig:1} A typical configuration}
 \end{figure}

We are considering the situation: $D$ is placed {\it behind} $D_0$ in the sence that one cannot see 
$D$ directly from $B$ because of the existence of $D_0$.
See Figure \ref{fig:1} for a typical configuration of $B$, $D_0$ and $D$ in two-space dimensions.
The $B$ is a model of the place where the wave is generated and observed.

In this paper, the governing equations of the wave is the following.

Let $0<T<\infty$.  Consider the following initial boundary value problem for the classical wave equation:
$$
\left\{
\begin{array}{ll}
\partial_t^2 u-\Delta u=0 & \text{in $(\Bbb R^3\setminus(\overline{D_0}\cup\overline D))\times\,]0,\,T[$,}\\
\\
\displaystyle
u(x,0)=0 & \text{in $\Bbb R^3\setminus(\overline{D_0}\cup\overline D)$,}\\
\\
\displaystyle
\partial_tu(x,0)=f(x) & \text{in $\Bbb R^3\setminus(\overline{D_0}\cup\overline D)$,}\\
\\
\displaystyle
\frac{\partial u}{\partial\nu}=0 & \text{on $\partial D_0\times\,]0,\,T[$,}\\
\\
\displaystyle
\frac{\partial u}{\partial\nu}=0 & \text{on $\partial D\times\,]0,\,T[$,}
\end{array}
\right.
\tag {1.1}
$$
where $f\in L^2(\Bbb R^3)$ with $\text{supp}\,f\subset\overline B$ and
$\nu$ denotes the outward normal to both $D_0$ on $\partial D_0$ and $D$ on $\partial D$.
We specify the meaning of the solution class of (1.1).
In short, it is the weak solution in \cite{DL} as considered in \cite{IW}.
More precisely, by Theorem 1 on p.558 in \cite{DL}, we know that there exists a unique $u=u(t)(\,\cdot\,)$ satisfying
$$\left\{
\begin{array}{l}
\displaystyle
u\in L^2(0,\,T; H^1(\Bbb R^3\setminus(\overline{D_0}\cup\overline D))),\\
\\
\displaystyle
u'\in L^2(0,\,T; H^1(\Bbb R^3\setminus(\overline{D_0}\cup\overline D))),\\
\\
\displaystyle
u''\in L^2(0,\,T; (H^1(\Bbb R^3\setminus(\overline{D_0}\cup\overline D)))')
\end{array}
\right.
$$
such that, for all $\phi\in H^1(\Bbb R^3\setminus(\overline{D_0}\cup\overline D))$,
$$\begin{array}{ll}
\displaystyle
<u''(t),\phi>
+\int_{\Bbb R^3\setminus(\overline{D_0}\cup\overline D)}\nabla u(t)\cdot\nabla\phi\,dx=0
&
\text{a.e. $t\in\,]0,\,T[$,}
\end{array}
$$
and $u(0)=0$, $u'(0)=f$.  In this paper, we say this $u$ is the solution of (1.1) and write $u=u(x,t)$.

This paper is concerned with the following inverse problem.

$\quad$

{\bf\noindent Problem 1.}
Fix $T$ (to be determined later), $B$ and $f$ (to be specified later).
Assume that $D_0$ is {\it known} and $D$ {\it unknown}.
Extract information about the location and shape of $D$ from the back-scattering data $u(x,t)$ given at all
$(x,t)\in\,B\times\,]0,\,T[$.

$\quad$

In \cite{IW2}, we have already shown that one can extract the distance 
of $B$ to $D_0\cup D$ provided $\text{ess}\,\inf_{x\in B}\,f(x)>0$. 
However, this contains both information about $D_0$ and $D$. 
Clearly the result does not yield a solution to Problem 1 if
the obstacle $D$ is placed behind $D_0$ where one can not see $D$ directly by using visible ray from $B$.
However, before considering this problem, first of all we have to consider
the following fundamental problem which is the subject of this paper.

$\quad$

{\bf\noindent Problem 2.}
Fix $T$ (to be determined later), $B$ and $f$ (to be specified later).
Assume that $D_0$ is {\it known}.
Can one know the {\it existence} of $D$ by observing  $u(x,t)$ given at all
$(x,t)\in\,B\times\,]0,\,T[$ ?

$\quad$

This is a mathematical formulation of the problem for obtaining information about the existence
of an unknown obstacle {\it behind} a known obstacle from the observation 
of a wave at a place where the wave is generated.  Note that we are considering
the case when one can not see directly the unknown obstacle by using visible ray.
In other words, this is a problem of finding a difference between two cases $D=\emptyset$ and $D\not=\emptyset$
in the wave $u(x,t)$ observed on $B\times\,]0,\,T[$.  
We are mainly thinking about the case when $D$ is in a so-called shadow region of a {\it sound-hard} obstacle $D_0$ from $B$
and seeking an {\it exact} and {\it constructive} criterion for the existence.
To the best of the author's knowledge there is no solution for Problems 1 and 2.

There are a lot of studies for clarifying the behaviour of a wave field in the shadow region which 
is a solution of the Helmholtz equation with {\it high frequency} in the exterior of a sound-hard obstacle
with a single connected component.
See \cite{B} for the exterior of a {\it convex obstacle} in the two-space dimensions,
whose shape depends on the positions of the observer and source;
\cite{f} for the justification of the leading profile of
the asymptotic expression of a wave field in the exterior of an arbitrary plane convex domain derived in \cite{M}
in the case when the source and point of observation to be situated on the boundary;
\cite{P} for an upper estimate of a wave field as frequency tends to infinity 
in the exterior of a {\it non-trapping} obstacle with a real analytic boundary
in $n$-space dimensions ($n\ge 2$).  See also the references in \cite{B,f,P} for related other works.
Those studies say that the wave field in the shadow region of a sound-hard obstacle from the source point is very weak 
as the frequency tends to infinity.

It is not clear that such studies yield the desired solution since: we observe the wave at the place where the wave is generated;
we have to consider the wave field outside the obstacle $D_0\cup D$ with $\overline{D_0}\cap\overline D=\emptyset$
which is {\it non convex} and may not non-trapping in general.
However, we expect that there should be an {\it impression}
of the existence of the unknown obstacle in the diffracted and reflected wave in time domain from $D$ since
the wave in the time domain shall be consists of several frequencies.

It should be pointed out also that there are some studies from engineering point of view which support
this expectation.  In \cite{NOKO, NKOO}, the detection problem of an unknown object behind another object (occlusion as called by them)
has been considered.  The problem is raised in robotics for object tracking.
They made use of audible sound (3.2kHz) generated by a speaker 
as a wave and catches the reflected sound by microphones placed where the object is line-of-site because of the occlusion.
They found that it is possible to give some estimate of the distance to the object when the object and occlusion are consists of two plates or
balls and placed in a simple configuration.  However their study is purely experimental and never formulate 
the problem as an inverse problem for the wave equation like Problems 1 and 2.

In this paper, we present an approach based on the {\it time domain enclosure method} developed
in \cite{IW, IW2, IW3, IW4, IW5, IE4, Thermo} which is a time domain version of the original enclosure method introduced in \cite{IE3, IE2, IE}.

First we introduce the indicator function of the time domain enclosure method in this paper.

{\bf\noindent Definition 1.1.}
Define
$$\begin{array}{ll}
\displaystyle
I_B(\tau;T)=\int_Bf(w-v)dx, & \tau>0,
\end{array}
\tag {1.2}
$$
where
$$\begin{array}{ll}
\displaystyle
w(x,\tau)
=\int_0^{T}e^{-\tau t}u(x,t)dt, & x\in\,\Bbb R^3\setminus(\overline{D_0}\cup\overline D),
\end{array}
\tag {1.3}
$$
and $v\in H^1(\Bbb R^3\setminus\overline{D_0})$ is the unique weak solution of
of the following problem:
$$\left\{
\begin{array}{ll}
\displaystyle
(\Delta-\tau^2)v+f=0 & \text{in $\Bbb R^3\setminus\overline{D_0}$,}\\
\\
\displaystyle
\frac{\partial v}{\partial\nu}=0 & \text{on $\partial D_0$.}
\end{array}
\right.
\tag {1.4}
$$
Since we assume that $D_0$ is known, in principle, $v$ is known.
We call the function $I_B(\tau;T)$ of independent variable $\tau$ the {\it indicator function}.
Define also
$$\displaystyle
J_{B}(\tau;D)=\int_{D}(\vert\nabla v\vert^2+\tau^2\vert v\vert^2)dx.
\tag {1.5}
$$
Note that both $I_B(\tau;T)$ and $J_B(\tau;D)$ depend also on $f$, however, for simplicity
of description we omit putting the symbol $f$ in them.

A solution to Problems 2 is the following criterion to make a descision whether $D=\emptyset$ or
not.  Note that at the present time we do not specify the form of $f$.

\proclaim{\noindent Theorem 1.1.}
We have:

(i)  if  $D=\emptyset$, then for all $T>0$ we have
$$\displaystyle
\lim_{\tau\longrightarrow\infty}e^{\tau T}I_B(\tau;T)=0;
$$

(ii)  if $D\not=\emptyset$ and $T$ satisfies
$$\displaystyle
\lim_{\tau\longrightarrow\infty}\,e^{\tau T}J_B(\tau;D)=\infty,
\tag {1.6}
$$
then we have
$$\displaystyle
\lim_{\tau\longrightarrow\infty}\,e^{\tau T}I_B(\tau;T)=\infty.
\tag {1.7}
$$

\endproclaim
Note that if some $T=T_0$ satisfies (1.6), then (1.6) is satisfied with all $T>T_0$. 

The following formula gives us a direct relationship between $D$ and the indicator function as $\tau\longrightarrow\infty$.

\proclaim{\noindent Theorem 1.2.}
Assume that $\partial D$ is $C^3$.  Then, (1.6) and (1.7) are equivalent and
we have
$$\displaystyle
\lim_{\tau\longrightarrow\infty}\frac{I_B(\tau;T)}{2J_B(\tau;D)}=1.
\tag {1.8}
$$

\endproclaim

At a first glance one thinks that, in order to clarify the asymptotic behaviour
of the indicator function $I_B(\tau;T)$ as $\tau\longrightarrow\infty$ one needs
to study the asymptotic behaviour on $B$ of the reduced wave $w$ given by (1.3) and $v$ 
as $\tau\longrightarrow\infty$.
However, Theorem 1.2 reduces the study to that of $v$ on $D$ as $\tau\longrightarrow\infty$.
Since $v$ is a solution of (1.4) and (1.4) is independent of $D$, this means that $D$ in (1.1) 
is eliminated in consideration. By formula (1.8), the leading profile of indicator function $I_B(\tau;T)$ as $\tau\longrightarrow\infty$ is completely
governed by that of $J_B(\tau;D)$ provided (1.6) or (1.7) is satified and $\partial D$
is sufficiently smooth.

Theorems 1.1 and 1.2 raise the following two questions:

$\bullet$  the existence of $T$ satisfying (1.6).

$\bullet$  if the existence has been ensured, then characterize the infimum of the set
$$\displaystyle
\left\{\frac{T}{2}\,\,\vert\,T>0,\, \lim_{\tau\longrightarrow\infty}\,e^{\tau T} J_B(\tau;D)=\infty\right\}.
$$

Here we give partial answers to these questions.
Before descrbing the answers, define
$$\displaystyle
\text{dist}(D,B)=\inf_{x\in B,\,y\in D}\,\vert x-y\vert.
$$
In what follows, from a technical reason, we specify the form of $f$:
$$
\displaystyle
f(x)=
\left\{
\begin{array}{ll}
\displaystyle
(\eta-\vert x-p\vert)^2g(x) & \text{if $x\in B$,}
\\
\\
\displaystyle
0 & \text{if $\Bbb R^3\setminus B$,}
\end{array}
\right.
$$
where $g\in C^2(\overline B)$ and satisfies $\inf_{x\in B}g(x)>0$; $p$ and $\eta$ are the center and radius
of $B$, respectively.  Note that $f$ belongs to $H^2(\Bbb R^3)$ and $\text{supp}\,f=\overline B$.

\proclaim{\noindent Theorem 1.3.}  
Assume that there exists an open ball $W$ such that $\overline{D_0}\subset W$ and
$$\displaystyle
B\cup D\subset\Bbb R^3\setminus\overline W.
\tag {1.9}
$$
Then, $J_B(\tau;D)$ satisfies (1.6) for all $T>2\,C\,\text{dist}\,(D,B)$ with
$$\displaystyle
C=\sqrt{2}\,\sqrt{\left(\frac{\pi}{4}\right)^2+1}.
\tag {1.10}
$$
Moreover, we have
$$\displaystyle
\liminf_{\tau\longrightarrow\infty}
\frac{1}{2\tau}
\log J_B(\tau;D)
\ge -C\,\text{dist}\,(D,B).
\tag {1.11}
$$

\endproclaim

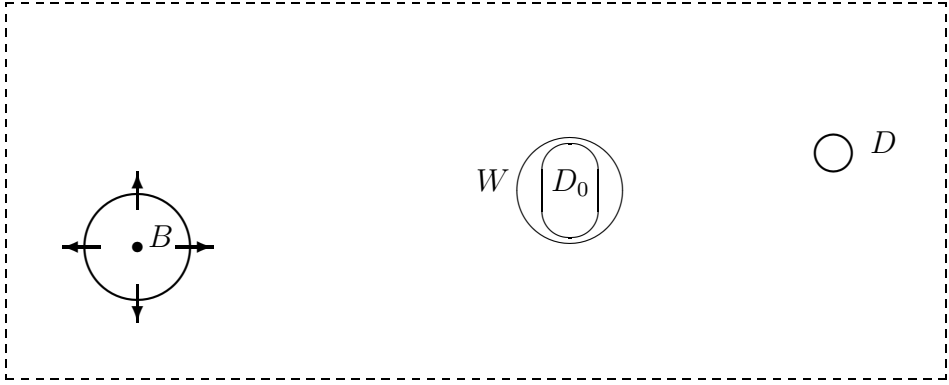
\begin{figure}
\setlength\unitlength{0.5truecm}
  \begin{picture}(10,10)(-3,-1)
  \put(0,0){\dashbox{0.2}(25,10){}}
  \put(15,5){\oval(1.5,2.5)}
  \put(12.5,5){$W$}
  \put(15,5){\circle{4}}
  \put(14.5,5){$D_0$}
  \put(3.5,3.5){\circle*{0.3}{$B$}}
  \thicklines
  \put(22,6){\circle{1}{$D$}}
  \put(3.5,3.5){\circle{3}}
  \put(4.5,3.5){\vector(1,0){1}}
  \put(3.5,4.5){\vector(0,1){1}}
  \put(3.5,2.5){\vector(0,-1){1}}
  \put(2.5,3.5){\vector(-1,0){1}}
  \end{picture}
 \caption{\label{fig:2} An illustration of ball $W$}
 \end{figure}

See Figure \ref{fig:2} for an illustration of the choice of ball $W$.
It should be emphasized that constant $C$ given by (1.10) is {\it independent} of $W$.
There is no assumption on the shape of $D_0$ together with its curvatures.
The $D_0$ can be the union of finitely many connected open subsets of $\Bbb R^3$ with $C^2$-boundaries
whose closures are mutually disjoint.

Needless to say we have the {\it explicit} upper bound 
$$\displaystyle
C\,\text{dist}\,(D,B)\ge
\inf
\left\{\frac{T}{2}\,\,\vert\,T>0,\, \lim_{\tau\longrightarrow\infty}\,e^{\tau T} J_B(\tau;D)=\infty\right\}
$$
provided $D$ satisfies (1.9) for an open ball $W$ satisfying $\overline{D_0}\subset W$.

{\bf\noindent Remark 1.1.}
Note that, without assuming the existence of $W$ in Theorem 1.3, we have for a positive constant $C'$
$$\displaystyle
\limsup_{\tau\longrightarrow\infty}\frac{1}{2\tau}\,\log J_B(\tau;D)\le -C'\text{dist}(D,B).
\tag {1.12}
$$
However, the dependence of $C'$ on $D_0$ is not clear.  See Appendix for the proof.
Anyway this ensures that the left-hand side on (1.11) is finite and negative.

As a corollary of Theorems 1.1, 1.3 and Proposition 2.1 together with Lemma 2.3 in Section 2 (for $\partial D\in C^2$ ) or
Theorems 1.1, 1.2 and 1.3 (for $\partial D\in C^3$) with the help of (1.12) we obtain

\proclaim{\noindent Corollary 1.1.}
Assume that $D$ satify (1.9) with an open ball $W$ satifying $\overline{D_0}\subset W$.
Let $C$ be the same constant given by (1.10).
Let $M>0$ satisfy $M>\text{dist}\,(D,B)$.
Fix $T\ge 2CM$.
Then, $D\not=\emptyset$ if and only if $\lim_{\tau\longrightarrow\infty}e^{\tau T}I_B(\tau;T)=\infty$.
Moreover, we have
$$\displaystyle
0>\liminf_{\tau\longrightarrow\infty}
\frac{1}{2\tau}
\log I_B(\tau;T)
\ge -C\,\text{dist}\,(D,B).
\tag {1.13}
$$

\endproclaim

If $D_0$ is an open ball, then one can always choose the $W$ in Theorem 1.3 by expanding $D_0$ a little bit.
Thus, we have

\proclaim{\noindent Corollary 1.2.}
Let $D_0$ be an open ball. Then, all the conclusions in Theorem 1.3 and Corollary 1.1 are valid.
\endproclaim

Therefore, we can detect any unknown obstacle $D$ {\it behind} an open ball by using a single wave over the finite time interval
$]0,\,T[$ with an arbitrary fixed $T$ satisfying $T\ge 2\,C\,M$ under a-priori assumption
$\text{dist}\,(D,B)<M$.  Moreover, (1.13) gives us the lower estimate of $\text{dist}\,(D,B)$ from the observed data:
$$\displaystyle
\text{dist}\,(D,B)>-\frac{1}{C}\liminf_{\tau\longrightarrow\infty}
\frac{1}{2\tau}
\log I_B(\tau;T).
$$
Note that we have $\text{dist}\,(D,B)=d_{\partial D}(p)-\eta$, where $d_{\partial D}(p)=\inf_{y\in\partial D}\,\vert y-p\vert$.  
Thus the estimate above gives us a lower estimate for $d_{\partial D}(p)$.

However, there are some situtation that $\overline D\subset W\setminus\overline{D_0}$ for any choice of ball $W$ with
constraint $\overline{D_0}\subset W$ and $\overline B\subset\Bbb R^3\setminus\overline W$.  
For example, consider the case when $D_0$ is an ellipsoid. See Figure \ref{fig:2}.
The result below covers this case under some a-priori assumtion on the possible location of $D$ relative to $D_0$.

We assume that $D_0$ is convex and thus $\overline D_0$ is also convex.
Then, given $x\in\Bbb R^3\setminus\overline D$ there exists a {\it unique point} $q(x)\in\partial D_0$ 
such that $\vert x-q(x)\vert=d_{\partial D_0}(x)$. Given $\alpha\in\,]-1,\,0]$ and $y\in B$ define
$$\displaystyle
V_{\alpha}(B;D_0)=\cap_{y\in B}\,V_{\alpha}(y;D_0),
$$
where
$$\displaystyle
V_{\alpha}(y;D_0)
=\left\{x\in\Bbb R^3\setminus\overline {D_0}\,\vert\,\nu_{q(x)}\cdot\nu_{q(y)}\ge\alpha\right\}.
$$

\proclaim{\noindent Theorem 1.4.}  
Let $D_0$ be convex. Let $B$ satisfy $\overline B\cap \overline{D_0}=\emptyset$.
Le $-1<\alpha\le 0$ and assume that $D$ satisfies
$$\displaystyle
D\subset V_{\alpha}(B;D_0).
\tag {1.14}
$$
Then, $J_B(\tau;D)$ satisfies (1.6) for all $T>2\,C(\alpha)\,\text{dist}\,(D,B)$ with
$$\displaystyle
C(\alpha)=\frac{\sqrt{2}}{1+\alpha}.
\tag {1.15}
$$
Moreover, we have
$$\displaystyle
\liminf_{\tau\longrightarrow\infty}
\frac{1}{2\tau}
\log J_B(\tau;D)
\ge -C(\alpha)\,\text{dist}\,(D,B).
\tag {1.16}
$$

\endproclaim

\begin{figure}
\setlength\unitlength{0.5truecm}
  \begin{picture}(10,10)(-3,-1)
  \put(0,0){\dashbox{0.2}(25,10){}}
  \thicklines
  \put(15,5){\circle{3}{$D_0$}}
  \put(10.5,6){$\nu_{q(y)}$}
  \put(13.6,5){\vector(-1,0){2}}
  \put(3.5,5){\circle*{0.1}$y$}
  \put(16.4,5.4){\line(6,1){7}}
  \put(16.4,4.6){\line(6,-1){7}}
  \put(8,2){$V_{\alpha}(y;D_0)$}
  \end{picture}
 \caption{\label{fig:3} An illustration of $V_{\alpha}(y;D_0)$}
 \end{figure}
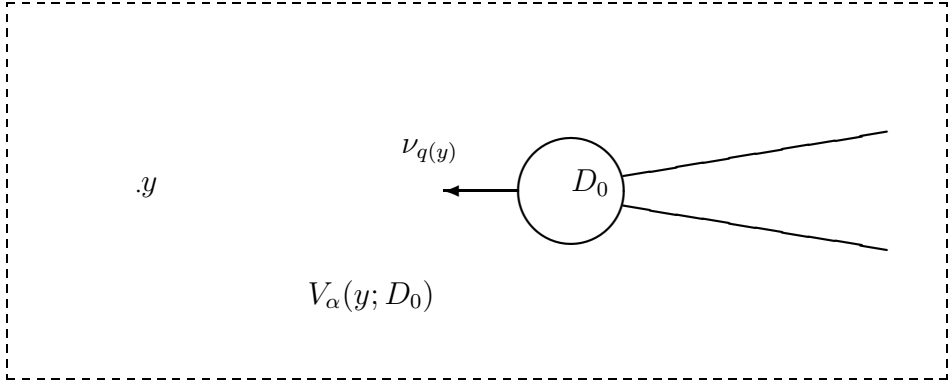

See also Figure \ref{fig:3} for an illustration of $V_{\alpha}(y;D_0)$ in the case when $D_0$ is an open ball.

\proclaim{\noindent Corollary 1.3.}
Let $D_0$ be convex.  Assume that $D$ satify (1.14) with $\alpha\in\,]-1, 0]$.
Let $C(\alpha)$ be the constant given by (1.15).
Let $M>0$ satisfy $M>\text{dist}\,(D,B)$.
Fix $T\ge 2C(\alpha)M$.
Then, $D\not=\emptyset$ if and only if $\lim_{\tau\longrightarrow\infty}e^{\tau T}I_B(\tau;T)=\infty$.
Moreover, we have
$$\displaystyle
0>\liminf_{\tau\longrightarrow\infty}
\frac{1}{2\tau}
\log I_B(\tau;T)
\ge -C(\alpha)\,\text{dist}\,(D,B).
\tag {1.17}
$$

\endproclaim

Thus, one can {\it hear} of the existence of an unknown obstacle behaind a known impenetrable {\it convex} obstacle
using a single wave over the finite time interval $]0,\,T[$ for an arbitrary {\it fixed} $T>2C(\alpha)M$
under the a-priori information (1.14) and $M>\text{dist}\,(D,B)$.
Moreover, (1.17) gives us a lower estimate of $\text{dist}\,(D,B)$ from the observed data.
However, if $D$ is placed in $\Bbb R^3\setminus V_{\alpha}(B;D)$ or a nearest point on $\overline D$ to $\overline B$ is
therein which is the case when $D$ is {\it close to the backside} of $D_0$ from $B$,
one can not say the behaviour of the indicator function at the present time.

Here we give a remark for general $D_0$.
Given $\epsilon>0$ define $(D_0)_{\epsilon}=\{x\in\Bbb R^3\setminus\overline{D_0}\,\vert\,d_{\partial D_0}(x)\ge \epsilon\}$.
Let $d_{\epsilon}(x,y)$ denote the infimum of the length of arcs in $(D_0)_{\epsilon}$ 
with endpoints $x$ and $y$.  If there is no arc in $(D_0)_{\epsilon}$ with enpoints $x$ and $y$, we define $d_{\epsilon}(x,y)=+\infty$.

If $D_0$ is not necessary convex, using the similar argument for the proof of Theorem 1.3, we have

\proclaim{\noindent Theorem 1.5.}
Assume that there exists a positive number $\epsilon$ such that
$B\cup D\subset (D_0)_{\epsilon}$ and
$$\displaystyle
\sup_{x\in B, y\in D}\,d_{\epsilon}(x,y)<\infty.
$$
Then, (1.6) is satisfied for all $T>2\sup_{x\in B,\,y\in D}\,d_{\epsilon}(x,y)$.

\endproclaim

Note that from this one gets a general result, however, which is {\it too rough} compared to Corollaries 1.1 and 1.3.

\proclaim{\noindent Corollary 1.4.}
Under the assumptions in Theorem 1.5 we have: $D\not=\emptyset$ if and only if $\lim_{\tau\longrightarrow\infty}e^{\tau T}I_B(\tau;T)=\infty$
for all $T>2\sup_{x\in B,\,y\in D}\,d_{\epsilon}(x,y)$.
Moreover, we have
$$\displaystyle
0>\liminf_{\tau\longrightarrow\infty}
\frac{1}{2\tau}
\log I_B(\tau;T)
\ge -2\sup_{x\in B,\,y\in D}\,d_{\epsilon}(x,y).
$$

\endproclaim
This is a solution of Problem 2 for $D_0$ without assuming the convexity nor a priori assumption (1.9)
/(1.14).  However, the intrinsic geometrical quantity $\sup_{x\in B\,,y\in D}d_{\epsilon}(x,y)$ is not useful 
for the {\it ranging} of unknown $D$
compared with the {\it Euclidean distance} $\text{dist}\,(D,B)$.

Finally we point out that the $v$ in (1.2) can be replaced with the $V$ given by
$$\displaystyle
V(x,\tau)=\int_0^T e^{-\tau t}u_0(x,t)dt,
$$
where $u_0$ is the weak solution of 
$$
\left\{
\begin{array}{ll}
\partial_t^2 u_0-\Delta u_0=0 & \text{in $\Bbb R^3\setminus\overline{D_0}\times\,]0,\,T[$,}\\
\\
\displaystyle
u_0(x,0)=0 & \text{in $\Bbb R^3\setminus\overline{D_0}$,}\\
\\
\displaystyle
\partial_tu_0(x,0)=f(x) & \text{in $\Bbb R^3\setminus\overline{D_0}$,}\\
\\
\displaystyle
\frac{\partial u_0}{\partial\nu}=0 & \text{on $\partial D_0\times\,]0,\,T[$.}
\end{array}
\right.
$$
More precisely, define
$$\begin{array}{ll}
\displaystyle
I'_B(\tau;T)=\int_Bf(w-V)dx, & \tau>0.
\end{array}
$$
It is easy to see that we have, as $\tau\longrightarrow\infty$
$$\displaystyle
\int_Bf(v-V)dx=O(\tau^{-1}e^{-\tau T})
$$
and hence
$$\displaystyle
I'_B(\tau;T)=I_B(\tau;T)+O(\tau^{-1}e^{-\tau T}).
$$
Thus, we have
\proclaim{\noindent Corollary 1.5.}
Theorems 1.1, 1.2 and Corollaries 1.1, 1.2, 1.3 and 1.4 remain vaild if $I_B(\tau;T)$ is replaced with $I'_B(\tau;T)$.
\endproclaim

This corollary shall be useful for making a test of the performance of the time domain enclosure method
since the $V$ on $B$ also can be obtained by an experiment as done in \cite{NOKO}.

The results obtained in this paper are concerned with an inverse obstacle problem goverend by the classical wave equation. 
It would be interested to consider also the corresponding problem for the wave governed by the Maxwell system.
See \cite{IMax} for the time domain enclosure method for the Maxwell system in the case when $D=\emptyset$ and 
$D_0$ is unknown.

A brief outline of this paper is as follows.  Theorem 1. 1 is proved in Section 2.
The key point is the asymptotic decomposition formula of the indicator function as $\tau\longrightarrow\infty$
as described in Proposition 2.1.
In Section 3 Theorem 1.2 is proved by showing the asymptotic equivalence of the first and second terms in the asymptotic decomposition formula
as stated in Proposition 3.1.
The proof of Proposition 3.1 employs the Lax-Phillips reflection argument developed in \cite{LP} originally for the analysis
of the right-end point of the {\it scattering kernel}.
Theorems 1.3 and 1.4 are proved in Section 4.
Those are based on Lemmas 4.1 and 4.2 which are classical results for the Gaussian lower estimate of 
the {\it heat kernels} in Euclidean space due to van den Berg \cite{V, V2, V3};
Lemma 4.3 which is a corollary of a general comparison theorem of two semigroups due to Ouhabaz \cite{Oh2};
Lemmas 4.4 and 4.5 concerning with the geometry of the exterior of a ball and that of the closed convex set, respectively.
Those are listed in Section 4 and the proof of Lemma 4.4 is given.
In Appendix first the proof of Lemma 4.5 which empolys a part of an argument due to Gyrya and Saloff-Coste \cite{GS}
is given.  Next (1.12) is proved by using a result from \cite{Oh}.

In what follows, to avoid heavy notation, we simply write
$J_B(\tau;T)=J(\tau)$.

\section{Proof of Theorem 1.1}

The function $w$ given by (1.3)  satisfies
$$
\left\{\begin{array}{ll}
\displaystyle
(\Delta-\tau^2)w+f=e^{-\tau T}F_{\tau} & \text{in $\Bbb R^3\setminus(\overline{D_0}\cup\overline D)$,}\\
\\
\displaystyle
\frac{\partial w}{\partial\nu}=0 & \text{on $\partial D_0$,}\\
\\
\displaystyle
\frac{\partial w}{\partial\nu}=0 & \text{on $\partial D$,}
\end{array}
\right.
\tag{2.1}
$$
where
$$\displaystyle
F_{\tau}(x)=\partial_tu(x,T)+\tau u(x,T),\,\,x\in\,\Bbb R^3\setminus(\overline{D_0}\cup\overline D).
\tag {2.2}
$$

In what follows we write $\epsilon=w-v$.
The following lemma is a simple application of integration  by parts and now well understood in the enlosure method in the time domain.

\proclaim{\noindent Lemma 2.1.}
It holds that
$$\begin{array}{ll}
\displaystyle
\int_Bf(w-v)dx &
\displaystyle
=\int_D(\vert\nabla v\vert^2+\tau^2\vert v\vert^2)\,dx
+\int_{\Bbb R^3\setminus(\overline{D_0}\cup\overline D)}
(\vert\nabla\epsilon\vert^2+\tau^2\vert\epsilon\vert^2)dx\\
\\
\displaystyle
&
\,\,\,
\displaystyle
+e^{-\tau T}
\int_{\Bbb R^3\setminus(\overline{D_0}\cup\overline D)}F_{\tau}\epsilon dx
-e^{-\tau T}
\int_{\Bbb R^3\setminus(\overline{D_0}\cup\overline D)}F_{\tau}v dx.
\end{array}
\tag {2.3}
$$

\endproclaim

From (1.4) one can easily obtain the following estimate:
$$\displaystyle
\int_{\Bbb R^3\setminus\overline{D_0}}(\vert\nabla v\vert^2+\tau^2\vert v\vert^2)dx=O(\tau^{-2}).
\tag {2.4}
$$

\proclaim{\noindent Proposition 2.1.}
We have, as $\tau\longrightarrow\infty$
$$\displaystyle
\int_Bf(w-v)dx
=J(\tau)+E(\tau)+O(\tau^{-1}e^{-\tau T}),
\tag {2.5}
$$
where
$$\displaystyle
E(\tau)
=\int_{\Bbb R^3\setminus(\overline{D_0}\cup\overline D)}(\vert\nabla\epsilon\vert^2+\tau^2\vert\epsilon\vert^2)dx.
$$

\endproclaim

{\it\noindent Proof.}
Rewrite (2.3) as
$$\begin{array}{l}
\displaystyle
\,\,\,\,\,\,
\int_{\Bbb R^3\setminus(\overline{D_0}\cup\overline D)}
\vert\nabla\epsilon\vert^2 dx
+\tau^2\int_{\Bbb R^3\setminus(\overline{D_0}\cup\overline D)}
\left\vert\epsilon-\frac{f-e^{-\tau T}F_{\tau}}{2\tau^2}\right\vert^2dx
+J(\tau)\\
\\
\displaystyle
=\frac{1}{4\tau^2}
\int_{\Bbb R^3\setminus(\overline{D_0}\cup\overline D)}\vert f-e^{-\tau T}F_{\tau}\vert^2dx
+e^{-\tau T}\int_{\Bbb R^3\setminus(\overline{D_0}\cup\overline D)} F_{\tau}vdx.
\end{array}
\tag {2.6}
$$
This gives
$$\displaystyle
\int_{\Bbb R^3\setminus(\overline{D_0}\cup\overline D)}
\vert\nabla\epsilon\vert^2 dx
\le
\frac{1}{4\tau^2}
\int_{\Bbb R^3\setminus(\overline{D_0}\cup\overline D)}\vert f-e^{-\tau T}F_{\tau}\vert^2dx
+e^{-\tau T}\int_{\Bbb R^3\setminus(\overline{D_0}\cup\overline D)}\vert F_{\tau}v\vert dx.
\tag {2.7}
$$
Since from (2.2) and (2.4) we have $\Vert F_{\tau}\Vert_{L^2(\Bbb R^3\setminus(\overline{D_0}\cup\overline D)}=O(\tau)$
and $\Vert v\Vert_{\Bbb R^3\setminus(\overline{D_0}\cup\overline D)}=O(\tau^{-2})$, respectively,
it follows from (2.7) that
$\displaystyle \Vert\nabla\epsilon\Vert^2_{L^2(\Bbb R^3\setminus(\overline{D_0}\cup\overline D))}=O(\tau^{-2})$.
Similarry, it follows from (2.6) again we obtain 
$\displaystyle \Vert\epsilon\Vert^2_{L^2(\Bbb R^3\setminus(\overline{D_0}\cup\overline D))}=O(\tau^{-4})$.
Therefore applying these together with (2.4) to (2.3) we obtain (2.5).

\noindent
$\Box$

Using integration by parts, one gets

\proclaim{\noindent Lemma 2.2.}
We have
$$\displaystyle
E(\tau)
=\int_{\partial D}\frac{\partial v}{\partial\nu}\epsilon dS
-e^{-\tau T}\int_{\Bbb R^3\setminus (\overline{D_0}\cup\overline D)}F_{\tau}\epsilon dx.
\tag {2.8}
$$

\endproclaim

From the proof of Proposition 2.1 and (2.8), we have, as $\tau\longrightarrow\infty$
$$\displaystyle
E(\tau)
=\int_{\partial D}\frac{\partial v}{\partial\nu}\epsilon\,dS+O(e^{-\tau T}\tau^{-1}).
$$
In particular, if $D=\emptyset$, we have
$$\displaystyle
E(\tau)=O(e^{-\tau T}\tau^{-1})
\tag {2.9}
$$
and  $J_B(\tau;D)=0$.
Thus (2.5) becomes
$$\displaystyle
I_B(\tau;T)=O(e^{-\tau T}\tau^{-1}).
$$
Thus (i) is valid.  (ii) is a consequence of the non negativity of $E(\tau)$ and (2.5).

{\bf\noindent Remark 2.1.}
Note that (2.9) is a quite {\it rough} estimate.
If $D=\emptyset$, then it follows from (1.4) and (2.1) that $\epsilon=w-v$ satisfies
$$\left\{
\begin{array}{c}
\displaystyle
(\triangle-\tau^2)\epsilon=e^{-\tau T}F_{\tau}\,\,\text{in}\,\Bbb R^3\setminus\overline{D_0},\\
\\
\displaystyle
\frac{\partial\epsilon}{\partial\nu}=0\,\,\text{on}\,\partial D_0.
\end{array}
\right.
$$
From this, one can easily obtain $\Vert\epsilon\Vert_{L^2(\Bbb R^3\setminus\overline{D_0})}=O(\tau^{-1}e^{-\tau T})$.
Then, (2.8) in the case when $D=\emptyset$ yields
$$\displaystyle
E(\tau)=O(e^{-2\tau T}).
$$
However, the remainder term on (2.5) has bound $O(e^{-\tau T}\tau^{-1})$ and
so no affects the result.  Note that the proof of Theorem 1.1 is quite simple since
the main tool is just {\it integration by parts}.  We never make use of the trace theorem.

As a rough relationship between $E(\tau)$ and $J(\tau)$ we obtain

\proclaim{\noindent Lemma 2.3.}
Given $\eta>0$ there exist a positive constant $C_{\eta}$ such that,
for all $\tau\ge\eta$
$$\displaystyle
E(\tau)
\le C_{\eta}\left\{\tau^2J(\tau)
+e^{-2\tau T}\right\}.
\tag {2.10}
$$

\endproclaim

{\it\noindent Proof.}
Solve
$$\left\{
\begin{array}{ll}
\displaystyle
(\Delta-\tau^2)\epsilon_1=0 & \text{in $\Bbb R^3\setminus(\overline{D_0}\cup\overline D)$,}\\
\\
\displaystyle
\frac{\partial\epsilon_1}{\partial\nu}=-\frac{\partial v}{\partial\nu}
&
\text{on $\partial D$,}\\
\\
\displaystyle
\frac{\partial\epsilon_1}{\partial\nu}=0 & \text{on $\partial D_0$.}
\end{array}
\right.
$$
Set
$$
\displaystyle
\alpha
=\int_{\Bbb R^3\setminus(\overline{D_0}\cup\overline D)}
(\vert\nabla\epsilon_1\vert^2+\tau^2\vert\epsilon_1\vert^2)dx.
$$
By the trace theorem, one can find a $\tilde{\epsilon_1}\in H^1(D)$ such that $\tilde{\epsilon_1}=\epsilon_1$ on $\partial D$
and satisfies $\Vert\tilde{\epsilon_1}\Vert_{H^1(D)}
\le C\Vert \epsilon_1\vert_{\partial D}\Vert_{H^{1/2}(\partial D)}$, where $C=C(D)>0$.
Integration by parts and $v\in H^2(D)$ give
$$\begin{array}{ll}
\displaystyle
\alpha
&
\displaystyle
=\int_{\partial D}\frac{\partial v}{\partial\nu}\epsilon_1dS\\
\\
\displaystyle
&
\displaystyle
=\int_{\partial D}\frac{\partial v}{\partial\nu}\tilde{\epsilon}dS
\\
\\
\displaystyle
&
\displaystyle
=\int_D((\Delta v)\tilde{\epsilon_1}+\nabla v\cdot\nabla\tilde{\epsilon_1})dx\\
\\
\displaystyle
&
\displaystyle
=\int_D(\nabla v\cdot\nabla\tilde{\epsilon_1}+\tau^2 v\tilde{\epsilon_1})dx\\
\\
\displaystyle
&
\displaystyle
\le
(\Vert\nabla v\Vert_{L^2(D)}+\tau^2\Vert v\Vert_{L^2(D)})\Vert\tilde{\epsilon_1}\Vert_{H^1(D)}\\
\\
\displaystyle
&
\displaystyle
\le
C(\Vert\nabla v\Vert_{L^2(D)}+\tau^2\Vert v\Vert_{L^2(D)})\Vert \epsilon_1\vert_{\partial D}\Vert_{H^{1/2}(\partial D)}.
\end{array}
$$
Let $\eta>0$.  By the trace theorem we have
$$\displaystyle
\Vert \epsilon_1\vert_{\partial D}\Vert_{H^{1/2}(\partial D)}
\le C_{\eta}\sqrt{\alpha}, \forall\tau\ge\eta.
$$

We have
$$\displaystyle
\Vert v\Vert_{L^2(D)}\le\tau^{-1}\sqrt{J(\tau)}
$$
and
$$\displaystyle
\Vert\nabla v\Vert_{L^2(D)}\le\sqrt{J(\tau)}.
$$

From these we obtain
$$\displaystyle
\alpha
\le
C_{\eta}'(\tau+1)\sqrt{J(\tau)}\sqrt{\alpha}
$$
and thus
$$\displaystyle
\alpha
\le( C_\eta')^2(\tau+1)^2\,J(\tau).
\tag {2.11}
$$

Next solve
$$\left\{
\begin{array}{ll}
\displaystyle
(\Delta-\tau^2)\epsilon_r=e^{-\tau T}F_{\tau} & 
\text{in $\Bbb R^3\setminus(\overline{D_0}\cup\overline D)$,}\\
\\
\displaystyle
\frac{\partial\epsilon_r}{\partial\nu}=0 & \text{on $\partial D$,}\\
\\
\displaystyle
\frac{\partial\epsilon_r}{\partial\nu}=0 & \text{on $\partial D_0$.}
\end{array}
\right.
$$
We have
$$\displaystyle
\int_{\Bbb R^3\setminus(\overline{D_0}\cup\overline D)}
(\vert\nabla\epsilon_r\vert^2+\tau^2\vert\epsilon_r\vert^2)dx
\le e^{-2\tau T}\tau^{-2}\int_{\Bbb R^3\setminus(\overline{D_0}\cup\overline D)}\vert F_{\tau}\vert^2 dx
\le C e^{-2\tau T}.
\tag {2.12}
$$

Since $\epsilon=\epsilon_1+\epsilon_r$, from (2.11) and (2.12) we obtain (2.10).

\noindent
$\Box$

\section{Proof of Theorem 1.2}

It suffices to prove the following proposition since we have (ii) of Theorem 1.1 and (2.5) in Proposition 2.1.

\proclaim{\noindent Proposition 3.1.}
Assume that $\partial D$ is $C^3$.
If $T$ satisfies (1.7), then (1.6) is satisfied and
we have
$$\displaystyle
\lim_{\tau\longrightarrow\infty}\,\frac{E(\tau)}{J(\tau)}=1.
\tag {3.1}
$$

\endproclaim

{\it\noindent Proof.}
Since $v$ satisfies $(\Delta-\tau^2)v=0$ in $D$, integration by parts yields the expression
$$\displaystyle
J(\tau)=\int_{\partial D}\frac{\partial v}{\partial\nu}\,vdS.
$$
Then, the boundary condition for $\epsilon$ on $\partial D$ gives
$$\displaystyle
J(\tau)=-\int_{\partial D}\frac{\partial \epsilon}{\partial\nu}\,vdS.
$$
A combination of this and (2.8) gives
$$\displaystyle
E(\tau)-J(\tau)
=\int_{\partial D}
\left(\frac{\partial v}{\partial\nu}\epsilon+
\frac{\partial\epsilon}{\partial\nu}v\right)dS
-e^{-\tau T}\int_{\Bbb R^3\setminus (\overline{D_0}\cup\overline D)}F_{\tau}\epsilon dx.
\tag {3.2}
$$

Let $\tilde{v}\in H^2(\Bbb R^3\setminus (\overline{D_0}\cup\overline D))$ satisfy
$$
\left\{
\begin{array}{ll}
\displaystyle
\tilde{v}=v\,\,\text{on $\partial D$,}\\
\\
\displaystyle
\frac{\partial\tilde{v}}{\partial\nu}
=-\frac{\partial v}{\partial\nu} & \text{on $\partial D$}
\end{array}
\right.
\tag {3.3}
$$
and
$$
\left\{
\begin{array}{ll}
\displaystyle
\tilde{v}=0 & \text{on $\partial D_0$,}\\
\\
\displaystyle
\frac{\partial\tilde{v}}{\partial\nu}
=0 & \text{on $\partial D_0$.}
\end{array}
\right.
\tag {3.4}
$$
Then the first term of the right-hand side on (3.2) becomes
$$\displaystyle
\int_{\partial D}
\left(-\frac{\partial\tilde{v}}{\partial\nu}\epsilon+
\frac{\partial\epsilon}{\partial\nu}\tilde{v}\right)dS.
$$
Integration by parts together with (3.4) gives
$$\begin{array}{c}
\displaystyle
-\int_{\partial D}\frac{\partial\tilde{v}}{\partial\nu}\epsilon dS
=\int_{\Bbb R^3\setminus(\overline{D_0}\cup\overline D)}\nabla\tilde{v}\cdot\nabla\epsilon dx
+\int_{\Bbb R^3\setminus(\overline{D_0}\cup\overline D)}(\Delta\tilde{v})\epsilon dx
\end{array}
$$
and
$$\begin{array}{c}
\displaystyle
-\int_{\partial D}\frac{\partial\epsilon}{\partial\nu}\tilde{v}dS
=\int_{\Bbb R^3\setminus(\overline{D_0}\cup\overline D)}\nabla\epsilon\cdot\nabla\tilde{v} dx
+\int_{\Bbb R^3\setminus(\overline{D_0}\cup\overline D)}(\Delta\epsilon)\tilde{v} dx.
\end{array}
$$
This together with the governing equation of $\epsilon$ yields
$$\displaystyle
\int_{\partial D}
\left(-\frac{\partial\tilde{v}}{\partial\nu}\epsilon+
\frac{\partial\epsilon}{\partial\nu}\tilde{v}\right)dS
=\int_{\Bbb R^3\setminus(\overline{D_0}\cup\overline D)}(\Delta-\tau^2)\tilde{v}\epsilon dx
-e^{-\tau T}\int_{\Bbb R^3\setminus(\overline{D_0}\cup\overline D)}F_{\tau}\tilde{v}dx.
$$
Therefore (3.2) takes the form
$$\displaystyle
E(\tau)-J(\tau)
=
\int_{\Bbb R^3\setminus(\overline{D_0}\cup\overline D)}(\Delta-\tau^2)\tilde{v}\cdot\epsilon dx
-e^{-\tau T}\left(\int_{\Bbb R^3\setminus(\overline{D_0}\cup\overline D)}F_{\tau}(\epsilon+\tilde{v})dx\right).
\tag {3.5}
$$
This is the starting point of the proof of (3.1).

In what follows $C$ denotes a positive constant.
Define
$$\displaystyle
\tilde{v}(x)=\phi(x)v(x^r),\,\,x\in\Bbb R^3\setminus\,(\overline{D_0}\cup\overline D).
\tag {3.6}
$$
where $\phi=\phi_{\delta}\in C^{\infty}(\Bbb R^3)$  with $0<\delta<\delta_0$ for a sufficiently small $\delta_0>0$
such that $0\le\phi\le 1$; $\phi(x)=1$ if $d_{\partial D}(x)<\delta$ and $\phi(x)=0$ if $d_{\partial D}(x)>2\delta$;
$\vert\nabla\phi(x)\vert\le C\delta^{-1}$; $\vert\nabla^2\phi(x)\vert\le C\delta^{-2}$;
$x^r=2q(x)-x$ for $x\in\Bbb R^3\setminus D$ with $d_{\partial D}(x)<2\delta_0$ where $q(x)$ denotes the unique point
on $\partial D$ such that $d_{\partial D}(x)=\vert x-q(x)\vert$.  We assume that $\partial D$ is $C^3$.  This ensures
that $q(x)$ is $C^2$ for $x\in\Bbb R^3\setminus D$ with $d_{\partial D}(x)<2\delta_0$ (see \cite{GT}).

We see that $\tilde{v}$ given by (3.6) satisfies (3.3) and (3.4) and has a compact support.

We have
$$\begin{array}{l}
\displaystyle
\,\,\,\,\,\,
\Delta(\tilde{v})-\phi(x)(\Delta v)(x^r)\\
\\
\displaystyle
=\phi(x)\sum_{i,j}a_{ij}(x)(\partial_i\partial_j v)(x^r)
+\sum_{j}b_j(x)(\partial_jv)(x^r)
+(\Delta\phi)(x)v(x^r),
\end{array}
\tag {3.7}
$$
where $\phi(x)a_{ij}(x)\in C^1_0(\Bbb R^3)$ and satisfies
$\vert \phi(x)a_{ij}(x)\vert\le Cd_{\partial D}(x)$; each $b_j(x)$
has the form
$$\displaystyle
b_j(x)=\sum_{k}b_{jk}(x)\partial_k\phi(x)+d_j(x)\phi(x)
$$
with $b_{jk}(x)\partial_k\phi(x)\in C^1_0(\Bbb R^3)$ and $d_j(x)\phi(x)\in C^0_0(\Bbb R^3)$.

Since $\overline D\cap(\overline B\cup\overline D_0)=\emptyset$,
if necessary, choosing a small $\delta_0$ again one has
$(\Delta-\tau^2)v=0$ in $d_{\partial D}(x)<2\delta_0$.
This gives $\phi(x)(\Delta v)(x^r)=\tau^2\phi(x)v(x^r)=\tau^2\tilde{v}(x)$ for all $x\in\Bbb R^3\setminus(\overline{D_0}\cup\overline D)$.
Substituting this into (3.7), we obtain
$$\displaystyle
(\Delta-\tau^2)\tilde{v}
=\phi(x)\sum_{i,j}a_{ij}(x)(\partial_i\partial_j v)(x^r)
+\sum_{j}b_j(x)(\partial_jv)(x^r)
+(\Delta\phi)(x)v(x^r).
$$
We substitute this into the first term of the right-hand side on (3.5)
and as previously done in \cite{IW3, IW4} one has
$$\begin{array}{c}
\displaystyle
\left\vert\int_{\Bbb R^3\setminus(\overline{D_0}\cup\overline D)}(\Delta-\tau^2)\tilde{v}\cdot\epsilon dx\right\vert
\le
C(\delta
+\delta^{-1}\tau^{-1}
+\delta^{-2}\tau^{-2})(J(\tau)E(\tau))^{1/2}.
\end{array}
$$
Choosing $\delta=\tau^{-1/2}$, we see that this right-hand side becomes
$\displaystyle C\tau^{-1/2}(J(\tau)E(\tau))^{1/2}$.

Moreover, it is easy to see that
$$\displaystyle
e^{-\tau T}\left\vert\int_{\Bbb R^3\setminus(\overline{D_0}\cup\overline D)}F_{\tau}(\epsilon+\tilde{v})dx\right\vert
\le Ce^{-\tau T}(E(\tau)^{1/2}+J(\tau)^{1/2}).
$$
Therefore from these together with (3.5)
we obtain
$$\displaystyle
\vert
E(\tau)-J(\tau)
\vert
\le
C\left(\tau^{-1/2}(J(\tau)E(\tau))^{1/2}+e^{-\tau T}(E(\tau)^{1/2}+J(\tau)^{1/2})\right).
\tag {3.8}
$$
We have, for $\epsilon>0$
$$\begin{array}{c}
\tau^{-1/2}(J(\tau)E(\tau))^{1/2}+e^{-\tau T}(E(\tau)^{1/2}+J(\tau)^{1/2})
\\
\\
\displaystyle
\le
\tau^{-1/2}\frac{J(\tau)+E(\tau)}{2}
+\frac{\epsilon^{-1}e^{-2\tau T}+\epsilon E(\tau)}{2}
+\frac{e^{-2\tau T}+J(\tau)}{2}\\
\\
\displaystyle
=\frac{1}{2}(\tau^{-1/2}+1)J(\tau)+
\frac{1}{2}(\tau^{-1/2}+\epsilon)E(\tau)
+\frac{1}{2}(1+\epsilon^{-1})e^{-2\tau T}.
\end{array}
\tag {3.9}
$$
Since $E(\tau)=E(\tau)-J(\tau)+J(\tau)\le\vert E(\tau)-J(\tau)\vert+J(\tau)$,
it follows from (3.8) and (3.9) that
$$\begin{array}{c}
\displaystyle
\left\{1-\frac{C}{2}(\tau^{-1/2}+\epsilon)\right\}
E(\tau)
\le 
\left\{1+\frac{C}{2}(\tau^{-1/2}+1)\right\}J(\tau)
+\frac{C}{2}(1+\epsilon^{-1})e^{-2\tau T}.
\end{array}
$$
From this we have 
$$\displaystyle
E(\tau)\le C(J(\tau)+e^{-2\tau T}).
\tag {3.10}
$$
A combination of this and (2.5) gives
$$\displaystyle
e^{\tau T}I_B(\tau;T)\le (C+1)e^{\tau T}J(\tau)+O(\tau^{-1}).
$$
Thus, (1.7) implies (1.6).

Moreover, from (3.10) we have
$$\begin{array}{l}
\,\,\,\,\,\,
\displaystyle
\tau^{-1/2}(J(\tau)E(\tau))^{1/2}+e^{-\tau T}(E(\tau)^{1/2}+J(\tau)^{1/2})\\
\\
\displaystyle
\le
\tau^{-1/2}\frac{J(\tau)+E(\tau)}{2}+e^{-\tau T}\frac{1+C(J(\tau)+e^{-2\tau T})}{2}
+e^{-\tau T}\frac{1+J(\tau)}{2}\\
\\
\displaystyle
=O(\tau^{-1/2}J(\tau)+e^{-\tau T}).
\end{array}
$$
Then from (3.8) again we obtain
$$\displaystyle
\vert E(\tau)-J(\tau)\vert
\le C(\tau^{-1/2}J(\tau)+e^{-\tau T})=CJ(\tau)\left(\tau^{-1/2}+\frac{1}{J(\tau)e^{\tau T}}\right).
$$
Applying (1.6) to this right-hand side, we immediately obtain (3.1).

\noindent
$\Box$

\section{Proof of Theorems 1.3 and 1.4}

\subsection{Some lemmas on the heat kernels}

In the following three lemmas, we do not assume that $D_0$ is convex.

\proclaim{\noindent Lemma 4.1(Theorem 2 in \cite{V2}).}
Let $K_0=K_0(x,y;t)$ be the Dirichlet heat kernel for the domain $\Bbb R^3\setminus\overline D_0$.
For $(x,y)\in ((D_0)_{\epsilon})^2$, $t>0$ and $\epsilon>0$
we have
$$\displaystyle
K_0(x,y;t)\ge\exp\left(-\frac{d_{\epsilon}(x,y)^2}{4t}\right)K_{\epsilon}(0,0;t),
$$
where $K_{\epsilon}(x,y;t)$ denotes the Dirichlet heat kernel for the open ball with radius $\epsilon$ centered at the origin.
\endproclaim

\proclaim{\noindent Lemma 4.2(Lemma 9 in \cite{V}).}
We have, for all $t>0$
$$\displaystyle
K_{\epsilon}(0,0;t)\ge (4\pi t)^{-3/2}\,e^{-9\pi^2t/(4\epsilon^2)}. 
$$
\endproclaim

Let $K_1=K_1(x,y;t)$ be the Neumann heat kernel for the domain $\Bbb R^3\setminus\overline{D_0}$
e.g., see Theorem 6.10 in \cite{Oh}.
The following lemma is a direct consequence of Corollary 4.3, Proposition 4.4 and Theorem 4.21 in \cite{Oh} which goes back to \cite{Oh2}.

\proclaim{\noindent Lemma 4.3.}
Let $f\in L^2(\Bbb R^3\setminus\overline{D_0})$ satify $f(y)\ge 0$ a.e. $x\in\Bbb R^3\setminus\overline{D_0}$.
We have, for all $t>0$ and $x\in\Bbb R^3\setminus\overline{D_0}$
$$\displaystyle
\int_{\Bbb R^3\setminus\overline{D_0}}K_0(x,y;t)f(y)dy\le \int_{\Bbb R^3\setminus\overline{D_0}}K_1(x,y;t)f(y)dy.
$$

\endproclaim

\subsection{On Theorem 1.3}

In what follows, given an arc $\gamma$ we denote by $L(\gamma)$ its length.
For the proof of Theorem 1.3, we employ the following simple fact.
\proclaim{\noindent Lemma 4.4.}
Let $U$ be an open ball.
Then, given $(x,y)\in (\Bbb R^3\setminus\overline U)^2$
there exists an arc $\gamma$ in $\Bbb R^3\setminus\overline U$ with endpoints $x$ and $y$ 
such that
$$\displaystyle
L(\gamma)
\le \sqrt{2}\,\sqrt{\left(\frac{\pi}{4}\right)^2+1}\,\vert x-y\vert.
\tag {4.1}
$$

\endproclaim

{\it\noindent Proof.}
Without loosing a generality, one may assume
that $d_{\partial U}(x)\le d_{\partial U}(y)$.  Let $\xi$ denote the center of $U$.
First consider the case when $(x-\xi)\cdot (y-\xi)<0$.
Let $x^*=2\xi-x_0$.   Let $U'$ denote the open ball centered at $\xi$ with radius $\vert x-\xi\vert$.
We have $\Bbb R^3\setminus U'\subset \Bbb R^3\setminus\overline U$ and $y\in\Bbb R^3\setminus U'$.
Let $y'$ denote the projection of $y$ onto the sphere $\partial U'$, that is
$$\displaystyle
y'=\xi+\vert x-\xi\vert\,\frac{y-\xi}{\vert y-\xi\vert}.
$$
Since $x$\,(North pole), $x^*$\,(South pole) and $y'$ are on the sphere $\partial U'$, one can find the {\it meridian} passing through
$y'$.  Let $\gamma_0$ denote tha part of the meridian that starts at $x$ and ends at $y'$.  
This is an arc in $\Bbb R^3\setminus\overline U$ with endpoints $x$ and $y'$.
It is clear that there is a unique point $x_0$ on $\gamma_0$ such that $(x_0-\xi)\cdot(x-\xi)=0$.
One gets
$$\displaystyle
L(\gamma_0')=\frac{\sqrt{2}\,\pi}{4}\,\vert x-x_0\vert,
\tag {4.2}
$$
where $\gamma_0'$ is the part of $\gamma_0$ that connects $x$ with $x_0$.

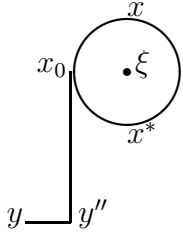
\begin{figure}
\setlength\unitlength{0.5truecm}
  \begin{picture}(12,12)(-10,-1)
  \thicklines
  \put(5,5){\circle{3}}
  \put(5,6.5){$x$}
  \put(5,3.1){$x^*$}
  \put(2.6,5){$x_0$}
  \put(3.7,1){$y''$}
  \put(1.8,1){$y$}
  \put(2.3,1){\line(1,0){1.2}}
  \put(3.5,1){\line(0,1){4}}
  \put(5,5){\circle*{0.2}$\xi$}
  \end{picture}
 \caption{\label{fig:4} An illustration of $x$, $x^*$, $y$, $x_0$ and $y''$}
 \end{figure}

Given two points $z_1$ and $z_2$ we denote by $[z_1,z_2]$ the line segment with endpoints $z_1$ and $z_2$.
One can find the unique arc $\gamma''=[x_0,y'']+[y'',y]$ such that $y''-x_0$ and $y-y''$ are parallel
to $x^*-\xi$ and $\xi-x_0$, respectively.  See Figure \ref{fig:4} for an irrustration of the choice.
Is is clear that $\gamma''$ is an arc in $\Bbb R^3\setminus\overline U$.
Since $\angle x_0y''y$ is right, we have
$$\displaystyle
L(\gamma'')=\vert x_0-y''\vert+\vert y''-y\vert\le\sqrt{2}\,\vert x_0-y\vert.
\tag {4.3}
$$
Now define $\gamma=\gamma_0'+\gamma''$. 
This is an arc in $\Bbb R^3\setminus\overline U$ with endpoints $x$ and $y$
and from (4.2) and (4.3) we have
$$\displaystyle
L(\gamma)
\le \sqrt{2}\,\left(\frac{\pi}{4}\vert x-x_0\vert+\vert x_0-y\vert\right).
\tag {4.4}
$$
We see that the triangle $\triangle xx_0x^*$ is a right triangle. 
Since the angle of two vectors $y-x_0$ and $x-x_0$
is greater than that of $x^*-x_0$ and $x-x_0$, one concludes that
the angle of two vectors $y-x_0$ and $x-x_0$
is greater than $\pi/2$.  This yields $(x-x_0)\cdot(y-x_0)<0$
and hence
$$\displaystyle
\sqrt{\vert x-x_0\vert^2+\vert x_0-y\vert^2}\le\vert x-y\vert.
$$
A combination of this and (4.4) yields (4.1).
If $(x-\xi)\cdot(y-\xi)\ge 0$, it is easy to find an arc $\gamma$ with endpoints $x$ and $y$ in $\Bbb R^3\setminus\overline U$,
lying in the half-space $\{z\in\Bbb R^3\,\vert\,(x-\xi)\cdot(z-\xi)\ge 0\}$ and having the length at most $\sqrt{2}\,\vert x-y\vert$.

\noindent
$\Box$

Now we give a proof of Theorem 1.3.
By (1.5), it suffices to prove that:  there exist real number $\kappa$, positive constants $C_1$ and $\tau_0$ such that,
for all $\tau\ge\tau_0$
$$\displaystyle
\tau^{\kappa}\Vert v\Vert_{L^2(D)}^2\ge C_1e^{-2\,C\,\text{dist}\,(D,B)\,\tau}.
\tag {4.5}
$$

The point is to introduce a solution of the initial bounday value problem for the heat equation
by using the semigroup theory \cite{Oh}.
The essence is the following.

Let $\cal {D}$ denotes the set of all functions $z\in H^1(\Bbb R^3\setminus\overline {D_0})$
satisfying that there exists a $\xi\in L^2(\Bbb R^3\setminus\overline{D_0})$ such that
$$\begin{array}{ll}
\displaystyle
\int_{\Bbb R^3\setminus\overline{D_0}}\,\nabla z\cdot\nabla\phi\,dx
=-\int_{\Bbb R^3\setminus\overline{D_0}}\,\xi\cdot\phi\,dx, & \forall\phi\in H^1(\Bbb R^3\setminus\overline{D_0}).
\end{array}
\tag {4.6}
$$
Note that $\xi=\Delta z$ in the sense of distribution.
Define 
$$\begin{array}{ll}
\displaystyle
Au=-\xi\in L^2(\Bbb R^3\setminus\overline{D_0}),& \forall u\in \cal{D}.
\end{array}
$$
It is known that $-A$ becomes the generator of a holomorphic semigroup $(e^{-tA})_{t\ge 0}$ on $L^2(\Bbb R^3\setminus\overline{D_0})$.
Given $\delta>0$ define $Z(x,t)=e^{-tA}f(x)$, $0\le t\le \delta$.  
Since $f\in\cal {D}$, we see that the map $Z(\,\cdot\,,t)\in\cal {D}$ for all $t\in\,[0,\,\delta]$
and the map:$[0,\,\delta]\longmapsto Z(\,\cdot\,,t)\in L^2(\Bbb R^3\setminus\overline{D_0})$ belongs to
$C^1([0,\,\delta],\,L^2(\Bbb R^3\setminus\overline{D_0}))$.
Set
$$\begin{array}{lll}
\displaystyle
v_{\delta}(x,,\tau)=\int_0^{\delta}e^{-\tau^2 t}Z(x,t)dt, & x\in\Bbb R^3\setminus\overline{D_0},
& \tau>0.
\end{array}
\tag {4.7}
$$
then $v(\,\cdot\,,\tau)$ belongs to $\cal{D}$.
Since $Z$ solves
$$\left\{
\begin{array}{ll}
\displaystyle
\frac{d}{dt}Z+AZ=0 & \text{in $[0,\delta]$,}
\\
\\
\displaystyle
Z(\,\cdot\,,0)=f,   &
\end{array}
\right.
$$
from (4.6) we have, in a weak sense
$$\left\{
\begin{array}{ll}
\displaystyle
(\Delta-\tau^2)v_{\delta}+f=e^{-\tau^2\delta}Z(x,\delta) & \text{in $\Bbb R^3\setminus\overline{D_0}$,}\\
\\
\displaystyle
\frac{\partial v_{\delta}}{\partial\nu}=0 & \text{on $\partial D_0$.}
\end{array}
\right.
$$
Set
$$\displaystyle
\epsilon_{\delta}=v-v_{\delta}.
$$
We have
$$\left\{
\begin{array}{ll}
\displaystyle
(\Delta-\tau^2)\epsilon_{\delta}=-e^{-\tau^2\delta}Z(x,\delta) & \text{in $\Bbb R^3\setminus\overline{D_0}$,}
\\
\\
\displaystyle
\frac{\partial\epsilon_{\delta}}{\partial\nu}=0 & \text{on $\partial D_0$.}
\end{array}
\right.
$$
Then, it is a routine to have the estimate
$$
\displaystyle
\int_{\Bbb R^3\setminus\overline{D_0}}
\left(\vert\nabla\epsilon_{\delta}\vert^2+\frac{\tau^2}{2}\vert\epsilon_{\delta}\vert^2\right)dx
\le
\frac{e^{-2\tau^2\delta}}{2\tau^2}
\int_{\Bbb R^3\setminus\overline{D_0}}\vert Z(x,\delta)\vert^2dx
$$
and hence
$$
\displaystyle
\Vert\nabla\epsilon_{\delta}\Vert^2_{L^2(\Bbb R^3\setminus\overline{D_0})}
+\tau^2\Vert\epsilon_{\delta}\Vert^2_{L^2(\Bbb R^3\setminus\overline{D_0})}
=O(\tau^{-2}e^{-2\tau^2\delta}).
$$
Since $v=v_{\delta}+\epsilon_{\delta}$, we obtain
$$\displaystyle
\Vert v\Vert_{L^2(D)}^2
\ge
\frac{1}{2}\Vert v_{\delta}\Vert_{L^2(D)}^2+O(\tau^{-4}e^{-2\tau^2\delta}).
\tag {4.8}
$$
Thus it suffices to give a lower estimate for $v_{\delta}$ given by (4.7).
For this we employ a Gaussian lower bound for the Dirichlet heat kernel in Lemmas 4.1 and 4.2 
combined with a relationship between two semigroups.

From the meaning of the heat kernel one has the expression
$$\displaystyle
Z(x,t)=\int_{\Bbb R^3\setminus\overline{D_0}}K_1(x,y;t)f(y)dy.
\tag {4.9}
$$
Then, by Lemma 4.3 together with $f(x)=0$ for $x\in\Bbb R^3\setminus B$, we have
$$\displaystyle
Z(x,t)\ge\int_BK_0(x,y;t)f(y)dy.
\tag {4.10}
$$
Let $0<t<\delta$ and $\epsilon>0$.   A combination of Lemmas 4.1 and 4.2 gives for $(x,y)\in((D_0)_{\epsilon})^2$,
$$\displaystyle
K_0(x,y;t)
\ge e^{-9\pi^2\delta/(4\epsilon^2)}
(4\pi t)^{-3/2}
\exp\left(-\frac{d_{\epsilon}(x,y)^2}{4t}\right).
\tag {4.11}
$$
Now let $\epsilon$ satisfy $\epsilon<\min(d,d_1)$ with $d=\text{dist}\,(B,D_0)$ and $d_1=\text{dist}\,(D,D_0)$.
Then, $B\cup D\subset (D_0)_{\epsilon}$ and hence (4.11) is valid for all $x\in D$ and $y\in B$.
Integrating both sides on (4.11) with repect to $y$ over $B$ and then using (4.10), we obtain,
for all $x\in D$
$$
\displaystyle
Z(x,t)\ge e^{-9\pi^2\delta/(4\epsilon^2)}
(4\pi t)^{-3/2}
\int_{B}
\exp\left(-\frac{d_{\epsilon}(x,y)^2}{4t}\right)f(y)dy.
$$
Again integrating both sides multiplied with $e^{-\tau^2 t}$ over time tinterval $]0,\,\delta[$, we obtain,
for all $x\in D$
$$\displaystyle
v_{\delta}(x,\tau)
\ge
e^{-9\pi^2\delta/(4\epsilon^2)}
\int_{B}
\left(\int_0^{\delta}(4\pi t)^{-3/2}
\exp\left(-\frac{d_{\epsilon}(x,y)^2}{4t}\right)
e^{-\tau^2 t}dt
\right)f(y)dy.
\tag {4.12}
$$
Now choose $\epsilon$ smaller in such a way that $\Bbb R^3\setminus\overline W\subset (D_0)_{\epsilon}$,
where $W$ is an open ball in the assumption of Theorem 1.3.
Let $x\in D$ and $y\in B$.  By (1.9) we have $(x,y)\in (\Bbb R^3\setminus\overline W)^2$.
Applying Lemma 4.4 to $U=W$, we can find an arc $\gamma$ with endpoints $x$ and $y$, lying in $\Bbb R^3\setminus\overline W$
and satisfies (4.1).  Since $\gamma$ lies also in $(D_0)_{\epsilon}$ we have
$$\displaystyle
d_{\epsilon}(x,y)\le C\vert x-y\vert,
\tag {4.13}
$$
where $C$ is the positive number given by (1.10).
Applying (4.13) to the right-hand side on (4.12), we obtain
$$\displaystyle
v_{\delta}(x,\tau)
\ge
e^{-9\pi^2\delta/(4\epsilon^2)}
\int_{B}
\left(\int_0^{\delta}(4\pi t)^{-3/2}
\exp\left(-C^2\frac{\vert x-y\vert^2}{4t}\right)
e^{-\tau^2 t}dt
\right)f(y)dy.
\tag {4.14}
$$
Since we have the formula
$$\begin{array}{ll}
\displaystyle
\int_0^{\infty}(4\pi t)^{-3/2}\exp\left\{\displaystyle -\frac{s^2}{4t}\right\}
e^{-\tau^2 t}dt
=\frac{e^{-\tau s}}{4\pi s}, & s>0,
\end{array}
\tag {4.15}
$$
we obtain, as $\tau\longrightarrow\infty$ uniformly with respect to $s>0$
$$\displaystyle
\int_0^{\delta}(4\pi t)^{-3/2}\exp\left\{\displaystyle -\frac{s^2}{4t}\right\}
e^{-\tau^2 t}dt
=\frac{e^{-\tau s}}{4\pi s}
+O(e^{-\tau^2\delta}).
$$
Applying this to the left-hand side on (4.14), we obtain
$$\displaystyle
v_{\delta}(x,\tau)\ge
e^{-9\pi^2\delta/(4\epsilon^2)}
\left(
\frac{1}{4\pi C}\int_{B}
\frac{e^{-C\tau\vert x-y\vert}}{\vert x-y\vert}f(y)dy
+O(e^{-\tau^2\delta})\right).
\tag {4.16}
$$
Here we make use of the special form of $f$:
$$\displaystyle
\int_{B}
\frac{e^{-C\tau\vert x-y\vert}}{\vert x-y\vert}f(y)dy
\ge C_1\int_{B}
\frac{e^{-C\tau\vert x-y\vert}}{\vert x-y\vert}(\eta-\vert y-p\vert)^2\,dy.
$$
A direct computation yields that, for sufficently large $\tau_0$ and all $\tau\ge\tau_0$
$$\begin{array}{ll}
\displaystyle
\int_{B}
\frac{e^{-C\tau\vert x-y\vert}}{\vert x-y\vert}(\eta-\vert y-p\vert)^2\,dy
\ge C_2\tau^{-3}\frac{e^{-C\tau\,(\vert x-p\vert-\eta)}}{\vert x-p\vert}, & \forall x\in\Bbb R^3\setminus\overline{B}.
\end{array}
$$
See, Lemma 2.3 in \cite{Thermo} and Appendix therein for a computation. 
And also one can show that for all $\tau>>1$,
$$\displaystyle
\tau^2 e^{2\tau\text{dist}\,(D,B)}
\int_De^{-2\tau(\vert x-p\vert-\eta)}dx\ge C_3.
$$
See, Lemma A.2 in \cite{Thermo}.  Note that the concreate value of the power of $\tau$ is not important.
A combination of this and (4.16) yields a lower estimate for $\Vert v_{\delta}\Vert_{L^2(D)}$ 
and finally (4.8) yields (4.5) with $\kappa=8$.

\subsection{On Theorem 1.4}

In the following lemma we make use of the assumption that $D_0$ is convex.

\proclaim{\noindent Lemma 4.5.}
Assume that $D_0$ is convex. Let $-1<\alpha\le 0$.
Then, for all $\epsilon>0$ and $(x,y)\in ((D_0)_{\epsilon})^2$ 
satisfying
$$\displaystyle
\nu_{q(x)}\cdot\nu_{q(y)}\ge\alpha
$$
we have
$$\displaystyle
d_{\epsilon}(x,y)\le C(\alpha)\,\vert x-y\vert.
\tag {4.17}
$$

\endproclaim

We do not think that the estimate (4.17) is optimal.  Say take $\alpha\downarrow -1$.
For the proof we employ a part of an argument done in the proof of Proposition 6.16 in \cite{GS}.
In the proposition they showed that the complement of an arbitrary closed convex set is {\it inner uniform} (see \cite{GS}).
The proof of Lemma 4.5 is given in Appendix.

Once we have Lemma 4.5, Theorem 1.4 can be proved along the same way as Theorem 1.3
by replacing $C$ in (4.13) with $C(\alpha)$ given by (1.15).
So we omit to describe the proof of Theorem 1.4.

$$\quad$$

\centerline{{\bf Acknowledgement}}

This research was partially supported by the Grant-in-Aid for
Scientific Research (C)(No. 17K05331) of Japan  Society for the
Promotion of Science.

\section{Appendix.}

\subsection{Proof of Lemma 4.5}
We apply an argument for the first part of the proof of Proposition 6.16 in \cite{GS}
done in the case when $\alpha=-1/\sqrt{2}$ to the case $\alpha$ general.

Let $(x,y)\in ((D_0)_{\epsilon})^2$.  This means that $(x,y)\in\,(\Bbb R^3\setminus\overline {D_0})^2$,
$d_{\partial D_0}(x)\ge\epsilon$ and $d_{\partial D_0}(y)\ge\epsilon$.
Let $H_x=\{z\in\Bbb R^3,\vert\,(z-x)\cdot\nu_{q(x)}\ge 0\}$ and $H_y=\{z\in\Bbb R^3,\vert\,(z-y)\cdot\nu_{q(y)}\ge 0\}$.
Since $D_0$ is convex, we see that $H_x\cup H_y\subset (D_0)_{\epsilon}$.

First we consider the case when $y\in H_x$ or $x\in H_y$.
Then, the straight line segment $[x,y]$ is contained in $H_x$ or $H_y$ and thus in $(D_0)_{\epsilon}$.
Thus we have 
$$
\displaystyle
d_{\epsilon}(x,y)=\vert x-y\vert.
\tag {A.1}
$$

Next consider the case when $y\notin H_x$ and $x\notin H_y$.  This means that
$(y-x)\cdot\nu_{q(x)}<0$ and $(y-x)\cdot\nu_{q(y)}>0$.  Thus it must have
$\nu_{q(x)}\cdot\nu_{q(y)}<1$.

Define
$$\left\{
\begin{array}{l}
\displaystyle
\mbox{\boldmath $t$}(x)=
\frac{1}
{\displaystyle\sqrt{1-(\nu_{q(x)}\cdot\nu_{q(y)})^2}}
\left(\nu_{q(y)}-(\nu_{q(y)}\cdot\nu_{q(x)})\nu_{q(x)}\right),\\
\\
\displaystyle
\mbox{\boldmath $t$}(y)=\frac{1}
{\displaystyle\sqrt{1-(\nu_{q(x)}\cdot\nu_{q(y)})^2}}
\left(\nu_{q(x)}-(\nu_{q(x)}\cdot\nu_{q(y)})\nu_{q(y)}\right).
\end{array}
\right.
$$
The unit vectors $\mbox{\boldmath $t$}(x)$ and $\mbox{\boldmath $t$}(y)$ satisfy
$$\displaystyle
\nu_{q(x)}\cdot\mbox{\boldmath $t$}(x)=\nu_{q(y)}\cdot\mbox{\boldmath $t$}(y)=0.
\tag {A.2}
$$
Moreover, we have
$$\displaystyle
\mbox{\boldmath $t$}(x)\cdot\mbox{\boldmath $t$}(y)=-\nu_{q(x)}\cdot\nu_{q(y)}.
\tag {A.3}
$$
Thus, $\mbox{\boldmath $t$}(x)$ and $\mbox{\boldmath $t$}(y)$ are linearly independent.
Therefore one can write
$$\displaystyle
y-x=a\mbox{\boldmath $t$}(x)+b\mbox{\boldmath $t$}(y)+c\mbox{\boldmath $t$}(x)\times\mbox{\boldmath $t$}(y),
\tag {A.4}
$$
where $a$, $b$, $c$ are real numbers depending on $x$ and $y$.
Moreover, we have $b<0$ and $a>0$ since $(y-x)\cdot\nu_{q(x)}<0$ and $(y-x)\cdot\nu_{q(y)}>0$.

From (A.3) and (A.4) one has
$$
\displaystyle
\vert y-x\vert^2=a^2+b^2-2ab\nu_{q(x)}\cdot\nu_{q(y)}+c^2.
\tag {A.5}
$$
Consder the case when $0\le\nu_{q(x)}\cdot\nu_{q(y)}$.  Since $ab<0$, from (A.5) we have
$$\displaystyle
\vert y-x\vert\ge \sqrt{a^2+b^2+c^2}.
\tag {A.6}
$$
Consider the case when $\alpha\le\nu_{q(x)}\cdot\nu_{q(y)}<0$.
Rewrite (A.5) as
$$\begin{array}{ll}
\displaystyle
\vert y-x\vert^2 & \displaystyle
=\left(\begin{array}{cc}
\displaystyle
1 & \displaystyle -\nu_{q(x)}\cdot\nu_{q(y)}\\
\\
\displaystyle
-\nu_{q(x)}\cdot\nu_{q(y)} &
\displaystyle
1
\end{array}
\right)
\left(\begin{array}{c}
a\\
\\
\displaystyle
b
\end{array}
\right)
\cdot
\left(\begin{array}{c}
a\\
\\
\displaystyle
b
\end{array}
\right)
+c^2.
\end{array}
\tag {A.7}
$$
Since we have
$$\displaystyle
\left(\begin{array}{cc}
\displaystyle
1 & \displaystyle -\nu_{q(x)}\cdot\nu_{q(y)}\\
\\
\displaystyle
-\nu_{q(x)}\cdot\nu_{q(y)} &
\displaystyle
1
\end{array}
\right)
\left(\begin{array}{c}
a\\
\\
\displaystyle
b
\end{array}
\right)
\cdot
\left(\begin{array}{c}
a\\
\\
\displaystyle
b
\end{array}
\right)
\ge
(1+\nu_{q(x)}\cdot\nu_{q(y)})(a^2+b^2),
$$
it follows from (A.7) that
$$\displaystyle
\vert x-y\vert\ge (1+\alpha)\sqrt{a^2+b^2+c^2}.
\tag {A.8}
$$

Now consider the arc $\gamma$ with end points $x$ and $y$ made of the three straight line segments
$$\displaystyle
\gamma=[x, x+a\mbox{\boldmath $t$}(x)]
+[x+a\mbox{\boldmath $t$}(x), y-b\mbox{\boldmath $t$}(y)]
+[y-b\mbox{\boldmath $t$}(y), y].
$$
From (A.6) and (A.8) we see that the length of $\gamma$ has the bound
$$\displaystyle
\vert a\vert+\vert b\vert+\vert c\vert
\le\sqrt{2}\sqrt{a^2+b^2+c^2}
\le
\frac{\sqrt{2}}{1+\alpha}\vert x-y\vert.
\tag {A.9}
$$

From (A.2) we see that the segments $[x, x+a\mbox{\boldmath $t$}(x)]$ and $[y-b\mbox{\boldmath $t$}(y), y]$
are contained in the planes $\pi_x=\partial H_x$ and $\pi_y=\partial H_y$, repectively. And also from (A.2) and (A.3) we have
$[x+a\mbox{\boldmath $t$}(x), y-b\mbox{\boldmath $t$}(y)]\subset\pi_x\cap\pi_y$.
Since $D_0$ is convex, one can conclude that $\gamma$ is contained in $(D_0)_{\epsilon}$.
Hence we have $d_{\epsilon}(x,y)\le \vert a\vert+\vert b\vert+\vert c\vert$.
this together with (A.1) and (A.9) gives (4.17).
This completes the proof of Lemma 4.5 with $C(\alpha)$ given by (1.15).

\subsection{Proof of (1.12)}

For $v_{\delta}$ given by (4.7), from (4.9) we have
$$\displaystyle
\nabla v_{\delta}(x,\tau)=\int_0^{\delta}e^{-\tau^2 t}\left(\int_B\nabla_xK_1(t,x,y)f(y)dy\right)dt.
\tag {A.10}
$$

Applying \cite{Oh}, (1) in Theorem 6.19 on page 185 to the present situation, we have
$$\begin{array}{lll}
\displaystyle
\int_{\Bbb R^3\setminus\overline{D_0}}\vert\nabla_xK_1(t,x,y)\vert^2 e^{\beta\vert x-y\vert^2/t}dx\le C_3t^{-1}, & t>0, & \text{a.e.}\,y\in\Bbb R^3\setminus\overline{D_0},
\end{array}
\tag {A.11}
$$
where $\beta$ and $C_3$ are positive constants.  Since $\overline D\subset\Bbb R^3\setminus\overline{D_0}$,
from (A.11) one has, for all $t>0$ and a.e. $y\in B$
$$\begin{array}{l}
\displaystyle
\,\,\,\,\,\,\int_D\vert\nabla_xK_1(t,x,y)\vert^2dx
\\
\\
\displaystyle
=
\int_D e^{-\beta\vert x-y\vert^2/t}\times
\vert\nabla_xK_1(t,x,y)\vert^2e^{\beta\vert x-y\vert^2/t}dx\\
\\
\displaystyle
\le
e^{-\beta\text{dist}\,(D,B)^2/t}
\int_D\vert\nabla_xK_1(t,x,y)\vert^2e^{\beta\vert x-y\vert^2/t}dx
\\
\\
\displaystyle
\le C_3e^{-\beta\text{dist}\,(D,B)\vert^2/t}t^{-1}.
\end{array}
\tag {A.12}
$$
(A.10) gives
$$\begin{array}{ll}
\displaystyle
\vert\nabla v_{\delta}(x,\tau)\vert
&
\displaystyle
\le\int_0^{\delta}e^{-\tau^2t}
\int_B\vert\nabla_xK_1(t,x,y)\vert f(y)dydt\\
\\
\displaystyle
&
\displaystyle
\le
\vert B\vert^{1/2}C_4
\int_0^{\delta}e^{-\tau^2 t}
\left(\int_B\vert\nabla _xK_1(t,x,y)\vert^2dy\right)^{1/2}dt\\
\\
\displaystyle
&
\le
\displaystyle
\vert B\vert^{1/2}\delta^{1/2}C_4
\left(\int_0^{\delta}e^{-2\tau^2 t}
\int_B\vert\nabla _xK_1(t,x,y)\vert^2 dydt\right)^{1/2}
\end{array}
$$
and thus
$$\begin{array}{l}
\displaystyle
\,\,\,\,\,\,
\int_D\vert\nabla v_{\delta}(x,\tau)\vert^2 dx
\\
\\
\displaystyle
\le\vert B\vert\delta C_4^2
\int_D dx
\int_0^{\delta}e^{-2\tau^2 t}
\int_B\vert\nabla _xK_1(t,x,y)\vert^2dydt\\
\\
\displaystyle
=\vert B\vert\delta C_4^2\int_0^{\delta}e^{-2\tau^2 t}dt
\int_Bdy
\int_D\vert\nabla_x K_1(t,x,y)\vert^2dx.
\end{array}
$$
A combination of this and (A.12) gives
$$\begin{array}{ll}
\displaystyle
\Vert\nabla v_{\delta}\Vert_{L^2(D)}^2
&
\displaystyle
\le \vert B\vert^2\delta C_4^2C_3\int_0^{\delta} e^{-2\tau^2 t}e^{-\beta\,\text{dist}\,(D,B)^2/t}t^{-1}dt\\
\\
\displaystyle
&
\displaystyle
\le\vert B\vert^2\delta^{3/2} C_4^2C_3\int_0^{\delta} e^{-2\tau^2 t}e^{-\beta\,\text{dist}\,(D,B)^2/t}t^{-3/2}dt.
\end{array}
$$
By (4.15) we have
$$\displaystyle
\int_0^{\infty}e^{-2\tau^2 t}e^{-\beta\text{dist}\,(D,B)^2/t}t^{-3/2}dt
=\sqrt{\frac{\pi}{\beta}}
\frac{e^{-2\sqrt{2}\sqrt{\beta}\,\text{dist}\,(D,B)\tau}}{\text{dist}(D,B)}.
$$
From these we obtain
$$\begin{array}{ll}
\displaystyle
\Vert\nabla v_{\delta}\Vert_{L^2(D)}^2
&
\displaystyle
\le\vert B\vert^2\delta^{3/2} C_4^2C_3
\left(\sqrt{\frac{\pi}{\beta}}
\frac{e^{-2\sqrt{2}\sqrt{\beta}\,\text{dist}\,(D,B)\tau}}{\text{dist}(D,B)}
+O(e^{-\tau^2\delta})\right)\\
\\
\displaystyle
&
\displaystyle
=O(e^{-2\sqrt{2}\sqrt{\beta}\,\text{dist}\,(D,B)\tau}+e^{-\tau^2\delta})
\end{array}
$$
and hence
$$\displaystyle
\Vert\nabla v_{\delta}\Vert_{L^2(D)}^2
=O(e^{-2\sqrt{2}\sqrt{\beta}\,\text{dist}\,(D,B)\tau}).
$$
Since $v=v_{\delta}+\epsilon_{\delta}$ and 
$\Vert\nabla\epsilon_{\delta}\Vert_{L^2(\Bbb R^3\setminus\overline{D_0})}^2=O(\tau^{-2}e^{-2\tau^2\delta})$,
we obtain
$$\displaystyle
\Vert\nabla v\Vert_{L^2(D)}^2
=O(e^{-2\sqrt{2}\sqrt{\beta}\,\text{dist}\,(D,B)\tau}).
\tag {A.13}
$$

Next we give an estimate for $\Vert v_{\delta}\Vert_{L^2(D)}^2$ from above.
By Theorem 6.10 in \cite{Oh} on page 171, we have
$$\displaystyle
\vert K_1(t,x,y)\vert
\le
C_5\,t^{-3/2}e^{-\alpha\vert x-y\vert^2/t}
$$
for all $t\in]\,0,\,\delta[$ and a.e. $(x,y)\in(\Bbb R^3\setminus\overline{D_0})^2$.  Here $\alpha$ and $C_5$ are positive constants.
Thus, we obtain, for a.e. $x\in D$ and a positive constatnt $C_6$
$$\begin{array}{ll}
\displaystyle
\vert v_{\delta}(x,\tau)\vert
&
\displaystyle
\le C_6\int_0^{\delta}e^{-\tau^2 t}\int_Be^{-\alpha\vert x-y\vert^2/t}dyt^{-3/2}dt\\
\\
\displaystyle
&
\displaystyle
\le
C_6\vert B\vert\int_0^{\delta}e^{-\tau^2 t}e^{-\alpha\,\text{dist}\,(D,B)^2/t}t^{-3/2}dt\\
\\
\displaystyle
&
\displaystyle
=C_6\vert B\vert
\left(
(4\pi)^{3/2}\frac{e^{-2\sqrt{\alpha}\,\text{dits}\,(D,B)\tau}}{4\pi\cdot 2\sqrt{\alpha}\,\text{dist}\,(D,B)}
+O(e^{-\tau^2\delta})
\right)
\\
\\
\displaystyle
&
\displaystyle
=O(e^{-2\sqrt{\alpha}\,\text{dist}(D,B)}).
\end{array}
$$
Therefore, we obtain
$$\displaystyle
\Vert v_{\delta}\Vert_{L^2(D)}^2=O(e^{-4\sqrt{\alpha}\,\text{dist}(D,B)\tau}).
$$
Since $\Vert\epsilon_{\delta}\Vert_{L^2(\Bbb R^3\setminus\overline{D_0)}}^2=O(\tau^{-4}e^{-2\tau^2\delta})$,
we finally obtain
$$\displaystyle
\Vert v\Vert_{L^2(D)}^2=O(e^{-4\sqrt{\alpha}\,\text{dist}(D,B)\tau}).
$$
A combination of this and (A.13) gives
$$\displaystyle
J(\tau)=O(e^{-2\gamma\,\text{dist}\,(D,B)\tau}),
$$
where $\gamma$ is a positive constant independent of $D$ and $B$.
From this we immediately obtain (1.12).

\vskip1cm
\noindent
e-mail address

ikehata@hiroshima-u.ac.jp


\begin{thebibliography}{99}





   
                  
\bibitem{B}  Bloom, C. O. and Matkowsky, B. J., On the validity of the geometrical theory of diffraction by convex cylinders,
             Arch. Rational Mech. Anal., {\bf 33}(1969), 71-90.
              
                                                                                
 
                
                   
                






\bibitem{DL} Dautray, R. and Lions, J-L.,
             Mathematical analysis and numerical methods for sciences and technology, 
             {\it Evolution Problems I}, vol. 5, Springer, Berlin, 1992.



\bibitem{f}  Filippov, V. B., 
             Rigorous justification of the shortwave asymptotic theory of diffraction 
             in the shadow zone,
             Zapiski Nauchnykh Seminarov Leningradskogo Otdeleniya Mathematicheskogo
             Instituta im. V. A.. Steklova Akademii Nauk SSSR,{\bf 34}(1973), 142-205.
             English Translation.





\bibitem{GT} Gilberg, D. and Trudinger, N. S.,
             Elliptic partial differential equations of second order, second. ed.,
             Springer-Verlag, Berlin, Heidelberg, New York, Tokyo, 1983.






\bibitem{GS} Gyrya, P. and Saloff-Coste, L.,
              Neumann and Dirichlet heat kernels in inner uniform domains,  Ast\'erisque, {\bf 336}, 2011.






\bibitem{IE3} Ikehata, M., Enclosing a polygonal cavity in a two-dimensional bounded domain from Cauchy data,
            Inverse Problems, {\bf 15}(1999), 1231-1241.




\bibitem{IE2} Ikehata, M., How to draw a picture of an unknown inclusion from boundary measurements.
            Two mathematical inversion algorithms,
            J. Inv. Ill-Posed Problems, {\bf 7}(1999), 255-271.




\bibitem{IE} Ikehata, M., Reconstruction of the support function for inclusion from boundary measurements,
            J. Inv. Ill-Posed Problems, {\bf 8}(2000), 367-378.







\bibitem{IW} Ikehata, M., The enclosure method for inverse obstacle scattering problems with dynamical data over
             a finite time interval, Inverse Problems, {\bf 26}(2010) 055010(20pp).





\bibitem{IW2} Ikehata, M., The enclosure method for inverse obstacle scattering problems with dynamical data over a finite time interval II.  Obstacles with a dissipative boundary or finite refractive index and back-scattering data, Inverse Problems, {\bf 28}(2012) 045010 (29pp).







\bibitem{IW3} Ikehata, M., An inverse acoustic scattering problem inside a cavity with dynamical back-scattering data, Inverse Problems, 
                           {\bf 28}(2012) 095016(24pp).







\bibitem{IW4} Ikehata, M., Extracting the geometry of an obstacle and a zeroth-order coefficient
of a boundary condition via the enclosure method using a single reflected wave over a finite time interval,
Inverse Problems, {\bf 30}(2014) 045011 (24pp).






\bibitem{IW5} Ikehata, M., On finding an obstacle embedded in the rough background medium via the enclosure method in the time domain,
                       Inverse Problems, {\bf 31}(2015) 085011 (21pp).



\bibitem{IMax}  Ikehata, M., The enclosure method for inverse obstacle scattering using a single electromagnetic wave in time domain,
                         Inverse Problems and Imaging, {\bf 10}(2016), 131-163.




                          
                         
                         
                         
\bibitem{IE4}  Ikehata, M., The enclosure method for inverse obstacle scattering over a finite time interval: IV.
                             Extraction from a single point on the graph of the response operator,
                             J. Inverse Ill-Posed Probl., 2017; aop, DOI: https://doi.org/10.1515/jiip-2016-0023

           
   
   
\bibitem{Thermo}  Ikehata, M., On finding a cavity in a thermoelastic body using a single displacement measurement over a finite time interval
                  on the surface of the body, arXiv:1706.02453v2 [math. AP].










\bibitem{LP}  Lax, P. D. and Phillips, R. S.,  The scattering of sound waves by an obstacle,
              Commun. Pure and Appl. Math., {\bf 30}(1977), 195-233.





\bibitem{M}  Molotkov, I.A., Field of a point source situated outside a convex curve,
             in: Problems of Mathematical Physics (in Russian), No.4, Izd. LGU (1970), 83-111.








\bibitem{NOKO} Niwa, H., Ogata, T., Komatani, K. and Okuno, H. G.,
               Detection and range finder of intercepted object by AH-based method using audible sound,
               IEICE Tech. Rep., 
               {\bf 106}(2006), no. 267, EA2006-48, 1-6, Sept. 2006, in Japanese.





\bibitem{NKOO} Niwa, H., Komatani, K., Ogata, T. and Okuno, H. G.,
               Distance estimation of hidden objects based on acoustical holography by applying acoustic diffraction of audible sound, 
               Proceedings of IEEE-RAS International Conference on Robotics and Automation (ICRA-2007), pp.423-428, (Apr. 2007). 
               doi:10.1109/ROBOT.2007.363823











\bibitem{Oh2} Ouhabaz, E. M., Invariance of closed convex sets and domination criteria for semigroups,
             Potential Anal., {\bf 5}(1996), 611-625.



\bibitem{Oh}  Ouhabaz, E. M.,  Analysis of heat equations on domains,
                    London Math. Soc. Monographs, vol. 31, Princeton University Press, 2004, Princeton and Oxford.




\bibitem{P}  Popov, G., 
             Some estimates of Green's functions in the shadow,
             Osaka J. Math., {\bf 24}(1987), 1-12.








\bibitem{V} van den Berg, M., Gaussian bounds for the Dirichlet heat kernel,
         J. Func. Anal., {\bf 88}(1990), 267-278.


\bibitem{V2} van den Berg, M., A Gaussian lower bound for the Dirichlet heat kernel,
         Bull. Lond. Math. Soc., {\bf 24}(1992), 475-477.



\bibitem{V3} van den Berg M., Subexponential behaviour of the Dirichlet heat kernel,
          J. Func. Anal., {\bf 198}(2003), 28-42.











\end{thebibliography}
\end{document}